\documentclass[a4paper]{amsart}
\usepackage{graphicx}
\usepackage{amsmath}
\usepackage{amssymb}
\usepackage{oldgerm}
\usepackage[all]{xy}

\newcommand{\Hom}{\operatorname{Hom}\nolimits}
\newcommand{\uHom}{\operatorname{\underline{Hom}}\nolimits}
\newcommand{\End}{\operatorname{End}\nolimits}

\renewcommand{\mod}{\operatorname{mod}\nolimits}
\newcommand{\stmod}{\operatorname{\underline{mod}}\nolimits}
\newcommand{\Ker}{\operatorname{Ker}\nolimits}

\newcommand{\Pd}{\operatorname{pd}\nolimits}
\newcommand{\gldim}{\operatorname{gldim}\nolimits}
\newcommand{\repdim}{\operatorname{repdim}\nolimits}

\newcommand{\Ext}{\operatorname{Ext}\nolimits}

\newcommand{\HH}{\operatorname{HH}\nolimits}

\newcommand{\rad}{\operatorname{rad}\nolimits}
\newcommand{\tr}{\operatorname{tr}\nolimits}
\newcommand{\res}{\operatorname{res}\nolimits}

\newcommand{\La}{\Lambda}
\newcommand{\Ga}{\Gamma}

\newcommand{\op}{\operatorname{op}\nolimits}

\newcommand{\add}{\operatorname{add}\nolimits}
\newcommand{\thick}{\operatorname{thick}\nolimits}

\newcommand{\Alt}{\mathcal{A}}
\newcommand{\cH}{\mathcal{H}}
\newcommand{\cS}{\mathcal{S}}

\newtheorem{theorem}{Theorem}[section]

\newtheorem{lemma}[theorem]{Lemma}
\newtheorem{proposition}[theorem]{Proposition}
\theoremstyle{definition}

\theoremstyle{definition}

\theoremstyle{definition}

\theoremstyle{definition}

\theoremstyle{definition}

\theoremstyle{definition}

\theoremstyle{remark}
\newtheorem*{remark}{Remark}
\theoremstyle{remark}
\newtheorem*{remarks}{Remarks}
\theoremstyle{definition}

\theoremstyle{definition}

\begin{document}
\title[Hecke algebras and symmetric groups]{The representation dimension of Hecke algebras and symmetric groups}
\author{Petter Andreas Bergh \& Karin Erdmann}
\address{Petter Andreas Bergh \newline Institutt for matematiske fag \\
  NTNU \\ N-7491 Trondheim \\ Norway}
\email{bergh@math.ntnu.no}
\address{Karin Erdmann \\ Mathematical Institute \\ 24-29 St.\ Giles \\ Oxford OX1 3LB \\ United Kingdom}
\email{erdmann@maths.ox.ac.uk}




\begin{abstract}
We establish a lower bound for the representation dimension
of all the classical Hecke algebras of types $A$, $B$ and $D$. For all the type $A$ algebras, and ``most" of the algebras of types $B$ and $D$, we also establish upper bounds. Moreover, we establish bounds for the representation dimension
of group algebras of some symmetric groups. 
\end{abstract}

\maketitle

\section{Introduction}\label{secintro}

The representation dimension of a finite dimensional algebra was
introduced by Auslander in \cite{Auslander1}. His aim was an
invariant which would somehow measure how far an algebra is from
having finite representation type. As a first step, he showed that
a non-semisimple algebra is of finite representation type if and
only if its representation dimension is exactly two, whereas it is
of infinite type if and only if the representation dimension is at
least three.

For more than three decades, no example was produced of an algebra
whose representation dimension exceeds three. However, in 2006
Rouquier showed in \cite{Rouquier2} that the representation
dimension of the exterior algebra on a $d$-dimensional vector
space is $d+1$, using the notion of the dimension of a
triangulated category (cf.\ \cite{Rouquier1}). Thus there do exist
finite dimensional algebras of arbitrarily large representation
dimension. Other examples illustrating this were subsequently
given in \cite{Bergh}, \cite{BerghOppermann1},
\cite{BerghOppermann2}, \cite{KrauseKussin}, \cite{Oppermann1},
\cite{Oppermann2}, \cite{Oppermann3}, \cite{OppermannMiemietz}.

Naturally, these papers focused on finding lower bounds for the representation dimension of various algebras. However, there does not exist a method for computing a good \emph{upper} bound. The best upper bound available so far was proved by Auslander himself: the representation dimension of a selfinjective algebra is at most its Loewy length. For some selfinjective algebras, this bound equals the representation dimension, but there also exist algebras for which the difference between this bound and the precise value is arbitrarily large. The simplest example is the selfinjective algebra $k[x]/(x^n)$ for $n \ge 2$. This has finite representation type, and its representation dimension is therefore $2$, whereas its Loewey length is $n$.

In this paper, we provide both an upper and a lower bound for the representation dimension of the Hecke
algebra $\mathcal{H}_q(A_{n-1})$ of type $A_{n-1}$, where $q$
is a primitive $\ell$th root of unity and the ground field is of
characteristic zero. In particular, we show that
$$[n/ \ell ] +1 \le \repdim \mathcal{H}_q(A_{n-1}) \le 2 [n/ \ell ]$$
whenever $\mathcal{H}_q(A_{n-1})$ is not semisimple (where $[r]$ denotes the integer part of a rational number $r$). These bounds are obtained by passing to a maximal $\ell$-parabolic subalgebra, which is just a tensor product of Brauer tree algebras and semisimple algebras. The proof carries over to some group algebras of symmetric groups, and consequently we obtain bounds also for such algebras. Namely, if $k$ is a perfect field of positive characteristic $p$, and $S_n$ is the $n$th symmetric group with $n<p^2$, then we show that the inequalities
$$[n/p] +1 \le \repdim kS_n \le 2 [n/p]$$
hold whenever $kS_n$ is not semisimple (i.e.\ when $n \ge p$). For these algebras, the bounds are obtained by passing to a Sylow $p$-subgroup of $S_n$: when $n<p^2$, this is an elementary abelian $p$-group of rank $[n/p]$.

We also establish lower bounds for all Hecke algebras of types $B$ and $D$, and these bounds are the same as for type $A$. Namely, if $\cH$ is either a Hecke algebra $\cH_{Q,q}(B_n)$ of type $B_n$, or a Hecke algebra $\cH_q(D_n)$ of type $D_n$, then we show that the inequality 
$$[n/ \ell ]+1 \le \repdim \cH $$
holds whenever $[n/ \ell ] \ge 1$. Moreover, when a certain polynomial expression in the parameters is nonzero (in the ground field), and $n$ is odd in type $D$, then we also provide an upper bound:
$$\repdim \cH \le 2[n/ \ell ].$$
As with the lower bound, this upper bound is the same as for type $A$.

\section{Comparing global dimensions}\label{seccomparing}

In this section we prove two results that compare the global dimensions of endomorphism rings of modules over two different algebras. They both apply to Hecke algebras of type $A$ and group algebras of symmetric groups. All modules are assumed to be finitely generated.

The first result considers the case of a subalgebra of an algebra.

\begin{theorem}\label{gldimsubalg}
Let $\La$ be a finite dimensional algebra, and suppose there exist
a subalgebra $\Ga$ and a $\Ga$-module $M$ such that:
\begin{enumerate}
\item the $\Ga$-module $\La \otimes_{\Ga}
M$ belongs to $\add_{\Ga} M$, \item the restriction map
$$\Ext_{\La}^1(X,Y) \xrightarrow{\res^{\La}_{\Ga}} \Ext_{\Ga}^1(X,Y)$$ is
injective for all $X,Y \in \add_{\La} ( \La \otimes_{\Ga} M)$.
\end{enumerate}
Then $\gldim \End_{\La} ( \La \otimes_{\Ga} M ) \le \gldim \End_{\Ga} ( M ).$
\end{theorem}

\begin{proof}
\sloppy Suppose the global dimension
of $\End_{\Ga} ( M )$ is finite, say $\gldim \End_{\Ga} ( M ) =d$.
Denote the induced $\La$-module $\La \otimes_{\Ga} M$ by
$M^{\La}$, and let $N_0$ be an indecomposable summand of
$M^{\La}$. Then $\Hom_{\La}( M^{\La},N_0 )$ is an indecomposable
projective $\End_{\La} ( M^{\La} )$-module, and every
indecomposable projective module for this endomorphism algebra is
of this form. Denote by $S(N_0)$ the simple $\End_{\La} ( M^{\La}
)$-module corresponding to $\Hom_{\La}( M^{\La},N_0 )$, and recall
that $\Hom_{\La}( M^{\La},-)$ induces an equivalence between
$\add_{\La} M^{\La}$ and the category of projective $\End_{\La}
(M^{\La} )$-modules. Therefore there exists an exact sequence
$$\cdots \to N_2 \to N_1 \to N_0$$
of $\La$-modules, in which each $N_i$ belongs to $\add_{\La} M^{\La}$, and such that
$$\cdots \to \Hom_{\La}( M^{\La},N_1) \to \Hom_{\La}( M^{\La},N_0) \to S(N_0) \to 0$$
is a projective resolution of $S(N_0)$. For each $i \ge 1$, denote
by $K_i$ the image of the map $N_i \to N_{i-1}$.

By restricting to $\Ga$ and using adjointness, we obtain an exact sequence
$$\cdots \to \Hom_{\Ga}( M,N_1) \to \Hom_{\Ga}( M,N_0) \to L \to 0$$
of $\End_{\Ga} ( M )$-modules. As $\Gamma$-modules, each $N_i$
belongs to $\add_{\Ga} M$, since $\add_{\Ga} M^{\La} \subseteq
\add_{\Ga} M$ by assumption. Therefore each $\Hom_{\Ga}( M,N_i)$
is a projective $\End_{\Ga} ( M )$-module, and the exact sequence
is a projective resolution of $L$. Since the projective dimension
of $L$ is at most $d$, the image of the map
$$\Hom_{\Ga}( M,N_d) \to \Hom_{\Ga}( M,N_{d-1}),$$
namely $\Hom_{\Ga}(M,K_d)$, is a projective $\End_{\Ga}(M)$-module. Consequently the short exact sequence
$$0 \to \Hom_{\Ga}(M,K_{d+1}) \to \Hom_{\Ga}( M,N_d) \to \Hom_{\Ga}( M,K_d) \to 0$$
of $\End_{\Ga}(M)$-modules splits. The functor $\Hom_{\Ga}( M,-)$
induces an equivalence between $\add_{\Ga} M$ and the category of
projective $\End_{\Ga} ( M )$-modules, hence the short exact
sequence
$$0 \to K_{d+1} \to N_d \to K_{d} \to 0$$
splits in $\mod \Ga$. By assumption, the restriction map
$$\Ext_{\La}^1(K_d, K_{d+1}) \xrightarrow{\res^{\La}_{\Ga}}
\Ext_{\Ga}^1(K_d, K_{d+1})$$ is injective, and so the short exact
sequence also splits in $\mod \La$. Then $K_d$ belongs to
$\add_{\La} M^{\La}$, hence in the projective resolution of the
$\End_{\La} ( M^{\La} )$-module $S(N_0)$ the image of the map
$$\Hom_{\La}( M^{\La},N_d) \to \Hom_{\La}( M^{\La},N_{d-1})$$
is projective. This shows that the projective dimension of every
simple $\End_{\La} ( M^{\La} )$-module is at most $d$, thus
proving the theorem.
\end{proof}

We shall also need the following result when we apply Theorem \ref{gldimsubalg} to Hecke algebras and group algebras. It gives sufficient conditions for an induced module to be a generator, that is, contain all the indecomposable projective modules as direct summands.

\begin{proposition}\label{addsubalg}
Let $\La$ be a finite dimensional algebra, and suppose there exist
a subalgebra $\Ga$ and a $\Ga$-module $M$ such that:
\begin{enumerate}
\item $\La$ is projective as a left $\Ga$-module, \item $M$ is a
generator in $\mod \Ga$.
\end{enumerate}
Then the $\La$-module $\La \otimes_{\Ga} M$ is a generator.
\end{proposition}

\begin{proof}
Since $\La$ is projective as a left $\Ga$-module, it belongs to
$\add_{\Ga} M$, and hence $\La \otimes_{\Ga} \La$ belongs to
$\add_{\La} ( \La \otimes_{\Ga} M)$. The surjective multiplication
map $\La \otimes_{\Ga} \La \xrightarrow{\mu} \La$ splits when
viewed as a map of left $\La$-modules, and so the left
$\La$-module $\La$ is a direct summand of $\La \otimes_{\Ga} \La$.
Therefore $\La$ must belong to $\add_{\La} ( \La \otimes_{\Ga}
M)$, and consequently $\La \otimes_{\Ga} M$ is a generator in
$\mod \La$.
\end{proof}

We turn now to the second result comparing global dimensions of endomorphism rings of modules over two different algebras. Whereas one of the algebras in Theorem \ref{gldimsubalg} was a subalgebra of the other, this is not necessarily the case in the following result. Still, the proofs are similar in nature.

Given two algebras $\La$ and $\Ga$, we say that $\La$ \emph{separably divides} $\Ga$ if there exist bimodules ${_{\La}X}_{\Ga}$ and ${_{\Ga}Y}_{\La}$, both projective on either side, such that the $\La$-$\La$-bimodule $\La$ is a direct summand of $X \otimes_{\Ga} Y$. If, in addition, the $\Ga$-$\Ga$-bimodule $\Ga$ is a direct summand of $Y \otimes_{\La} X$, then the algebras are \emph{separably equivalent} (cf.\ \cite{Linckelmann}). Obviously, if $\La$ and $\Ga$ are separably equivalent, then each of them separably divides the other. However, the converse does not seem to hold automatically: the bimodules involved need not be the same. Note that a group algebra is separably divided by the group algebra of any subgroup. Namely, if $k$ is a field, and $G$ is a group with a subgroup $H$, then the bimodules ${_{kG}X}_{kH}$ and ${_{kH}Y}_{kG}$ defined by $X = kG =Y$ do the trick: the $kH$-$kH$ bimodule $kH$ is a summand of $Y \otimes_{kG} X$. If, in addition, the subgroup $H$ is a Sylow subgroup, then $kG$ and $kH$ are separably equivalent, since in this case the $kG$-$kG$ bimodule $kG$ is a summand of $X \otimes_{kH} Y$. Similarly, a block algebra is separably equivalent to the group algebra of a defect group. 

\begin{theorem}\label{gldimseparably}
Let $\La$ and $\Ga$ be finite dimensional algebras, and suppose there exists
a $\Ga$-module $M$ such that:
\begin{enumerate}
\item $\La$ separably divides $\Ga$ through bimodules ${_{\La}X}_{\Ga}$ and ${_{\Ga}Y}_{\La}$,
\item $\Hom_{\La}(X, X \otimes_{\Ga} M) \in \add_{\Ga} M.$
\end{enumerate}
Then $\gldim \End_{\La} ( X \otimes_{\Ga} M ) \le \gldim \End_{\Ga} ( M ).$
\end{theorem}

\begin{proof}
Since $X$ is a projective left $\La$-module, the linear map
$$\Ext_{\La}^n(U,V) \xrightarrow{\Hom_{\La}(X,-)} \Ext_{\Ga}^n \left ( \Hom_{\La}(X,U) , \Hom_{\La}(X,V) \right )$$
is well-defined for all $n$ and all $\La$-modules $U,V$. Suppose an element $\eta \in \Ext_{\La}^n(U,V)$ maps to zero through this map, i.e.\ $\Hom_{\La}(X, \eta )=0$. Since $Y$ is projective as a left $\Ga$-module, the linear map
$$\Ext_{\Ga}^n(U',V') \xrightarrow{\Hom_{\Ga}(Y,-)} \Ext_{\La}^n \left ( \Hom_{\Ga}(Y,U') , \Hom_{\Ga}(Y,V') \right )$$
is well-defined for all $n$ and all $\Ga$-modules $U',V'$. Applying this map to the zero element $\Hom_{\La}(X, \eta )$, and using adjunction, we obtain
$$0 = \Hom_{\Ga}(Y, \Hom_{\La}(X, \eta ) ) \simeq \Hom_{\La}( X \otimes_{\Ga} Y, \eta).$$
Now since the $\La$-$\La$-bimodule $\La$ is a direct summand of $X \otimes_{\Ga} Y$, the extension $\eta$ is a direct summand of the extension $\Hom_{\La}( X \otimes_{\Ga} Y, \eta)$, hence $\eta =0$. Consequently, the linear map
$$\Ext_{\La}^n(U,V) \xrightarrow{\Hom_{\La}(X,-)} \Ext_{\Ga}^n \left ( \Hom_{\La}(X,U) , \Hom_{\La}(X,V) \right )$$
is injective for all $n$ and all $\La$-modules $U,V$.

Suppose the global dimension
of $\End_{\Ga} ( M )$ is finite, say $\gldim \End_{\Ga} ( M ) =d$.
Denote the induced $\La$-module $X \otimes_{\Ga} M$ by
$M^{\La}$, and let $S$ be a simple $\End_{\La} ( M^{\La} )$-module. The arguments used in the proof of Theorem \ref{gldimsubalg} show that there exists an exact sequence
$$\mathbb{S} \colon \cdots \to N_2 \to N_1 \to N_0$$
of $\La$-modules, with the following properties:
\begin{enumerate}
\item each $N_i$ belongs to $\add_{\La} M^{\La}$,
\item when applying $\Hom_{\La}( M^{\La},-)$ to $\mathbb{S}$, we obtain a projective resolution
$$\cdots \to \Hom_{\La}( M^{\La},N_1) \to \Hom_{\La}( M^{\La},N_0) \to S \to 0$$
of $S$ over $\End_{\La} ( M^{\La} )$.
\end{enumerate}
For each $i \ge 1$, denote by $K_i$ the image of the map $N_i \to N_{i-1}$.

Using adjointness, we obtain an isomorphism
\begin{eqnarray*}
\Hom_{\La}( M^{\La}, \mathbb{S} ) & = & \Hom_{\La}( X \otimes_{\Ga} M , \mathbb{S} ) \\
& \simeq & \Hom_{\Ga}(M, \Hom_{\La}( X, \mathbb{S} )).
\end{eqnarray*}
Consequently, the sequence $\Hom_{\La}( M^{\La}, \mathbb{S} )$ gives rise to an exact sequence
$$\cdots \to \Hom_{\Ga}(M, \Hom_{\La}( X, N_1)) \to \Hom_{\Ga}(M, \Hom_{\La}( X, N_0)) \to L \to 0$$
of $\End_{\Ga} ( M )$-modules. Since each $N_i$ belongs to $\add_{\La} M^{\La}$, and $\Hom_{\La}(X, M^{\La} ) \in \add_{\Ga} M$ by assumption, each $\End_{\Ga} ( M )$-module $\Hom_{\Ga}(M, \Hom_{\La}( X, N_i))$ is projective. Therefore
the sequence $\Hom_{\Ga}(M, \Hom_{\La}( X, \mathbb{S} ))$ is a projective resolution of $L$ as a module over $\End_{\Ga} ( M )$.

Since the global dimension of $\End_{\Ga} ( M )$ is $d$, the image of the map
$$\Hom_{\Ga}(M, \Hom_{\La}( X, N_d)) \to \Hom_{\Ga}(M, \Hom_{\La}( X, N_{d-1})),$$
namely $\Hom_{\Ga}(M, \Hom_{\La}( X, K_d))$, is a projective $\End_{\Ga} ( M )$-module. 
Therefore, when we apply $\Hom_{\Ga}(M,-)$ to the short exact sequence
\begin{equation*}\label{ES1}
0 \to \Hom_{\La}( X, K_{d+1}) \to \Hom_{\La}( X, N_d) \to \Hom_{\La}( X, K_d) \to 0 \tag{$\dagger$}
\end{equation*}
of $\Ga$-modules, the result is a split exact sequence of projective $\End_{\Ga} ( M )$-modules. But the functor $\mod \Ga \xrightarrow{\Hom_{\Ga}(M,-)} \mod \End_{\Ga} ( M )$ induces an equivalence between $\add_{\Ga} M$ and the category of projective $\End_{\Ga} ( M )$-modules, hence the sequence (\ref{ES1}) splits itself. It follows from the beginning of this proof that the linear map
$$\Ext_{\La}^1(K_d,K_{d+1}) \xrightarrow{\Hom_{\La}(X,-)} \Ext_{\Ga}^1 \left ( \Hom_{\La}(X,K_d) , \Hom_{\La}(X,K_{d+1}) \right )$$
is injective, and so the short exact sequence
$$0 \to K_{d+1} \to N_d \to K_d \to 0$$
splits in $\mod \La$. The module $K_d$ is then a summand of $N_d$, and since $N_d \in \add_{\La} M^{\La}$, we obtain $K_d \in \add_{\La} M^{\La}$. This implies that in the projective resolution of the $\End_{\La} ( M^{\La} )$-module $S$, 
the image $\Hom_{\La}( M^{\La},K_d)$ of the map
$$\Hom_{\La}( M^{\La},N_d) \to \Hom_{\La}( M^{\La},N_{d-1})$$
is projective, and consequently $\Pd_{\End_{\La} ( M^{\La} )} S \le d$. Since $S$ was an arbitrary simple module over $\End_{\La} ( M^{\La} )$, the global dimension of this endomorphism algebra is at most $d$, and the proof is complete.
\end{proof}

We end this section with the counterpart to Proposition \ref{addsubalg}

\begin{proposition}\label{addseparably}
Let $\La$ and $\Ga$ be finite dimensional algebras, and suppose there exists
a $\Ga$-module $M$ such that:
\begin{enumerate}
\item $\La$ separably divides $\Ga$ through bimodules ${_{\La}X}_{\Ga}$ and ${_{\Ga}Y}_{\La}$,
\item $M$ is a
generator in $\mod \Ga$.
\end{enumerate}
Then the $\La$-module $X \otimes_{\Ga} M$ is a generator.
\end{proposition}

\begin{proof}
Since $Y$ is projective as a left $\Ga$-module, it belongs to
$\add_{\Ga} M$, and hence the left $\La$-module $X \otimes_{\Ga} Y$ belongs to
$\add_{\La} ( X \otimes_{\Ga} M)$. But $\La$ is a direct summand of $X \otimes_{\Ga} Y$ as a bimodule, and in particular as a left $\La$-module. Therefore $\La$ belongs to $\add_{\La} ( X \otimes_{\Ga}
M)$.
\end{proof}

\section{Symmetric algebras}\label{secsymalg}

In this section, we record some properties of symmetric algebras.
All modules are assumed to be finitely generated left modules.
Recall that a finite dimensional algebra $\La$ over a field $k$ is
\emph{symmetric} if it is isomorphic as a bimodule to its $k$-dual
$D( \La ) = \Hom_k( \La, k)$. If this is the case, then fix such
an isomorphism $\La \xrightarrow{\phi} D( \La )$, and denote by $s
\in D( \La )$ the element $\phi (1)$. It is not hard to see that
\begin{eqnarray*}
s( xy ) &=&  s( yx ) \\
\phi ( x )( y ) &=& s( xy )
\end{eqnarray*}
for all $x,y \in \La$. The map $s$ is called a \emph{symmetrizing
form} for $\La$. Now suppose that $\Gamma$ is a \emph{parabolic}
subalgebra of $\La$ (cf.\ \cite{Broue}), that is, the following conditions are
satisfied:
\begin{enumerate}
\item $\Gamma$ is symmetric, \item the restriction of $s$ to
$\Gamma$ is a symmetrizing form for $\Gamma$, \item $\La$ is a
finitely generated projective left (equivalently, right)
$\Gamma$-module.
\end{enumerate}
Condition (2) is equivalent to the condition that, as a bimodule over $\Gamma$, the algebra $\La$ is a direct sum $\La = \Gamma \oplus B$, where $B \subseteq \Ker s$. In other words, the $\Gamma$-$\Gamma$-bimodule $\Gamma$ has a complement in $\La$ contained in the kernel of the symmetrizing form.

The multiplication map $\La \otimes_{\Ga} \La \xrightarrow{\mu}
\La$ and the bimodule isomorphism $\phi$ give rise to a
$\La$-$\La$-bimodule homomorphism $\La \to \La \otimes_{\Ga} \La$,
given as the composition
$$\La \xrightarrow{\phi} D( \La ) \xrightarrow{D( \mu )} D( \La
\otimes_{\Ga} \La ) \stackrel{\sim}{\longrightarrow} D( \La )
\otimes_{\Ga} D( \La ) \xrightarrow{\phi^{-1} \otimes \phi^{-1}}
\La \otimes_{\Ga} \La.$$ The image of the identity of $\La$ under
this homomorphism is the \emph{relative Casimir element}, and
denoted by $c_{\Ga}^{\La}$ (cf.\ \cite{Broue} and \cite{Linckelmann}). This element satisfies $x \cdot
c_{\Ga}^{\La} = c_{\Ga}^{\La} \cdot x$ for every $x \in \La$,
hence $\mu (c_{\Ga}^{\La})$ belongs to the center of $\La$. Write
$$c_{\Ga}^{\La} = \sum_{i=1}^n x_i \otimes y_i$$
for some elements $x_i,y_i \in \La$, and let $M$ and $N$ be two
$\La$-modules. Moreover, for a $\Ga$-homomorphism $M
\xrightarrow{f} N$, consider the $k$-linear map $M
\xrightarrow{\tr^{\La}_{\Ga} (f)} N$ defined by
$$\tr^{\La}_{\Ga}(f)(m) \stackrel{\text{def}}{=} \sum_{i=1}^n x_if(y_im)$$
for $m \in M$. It follows from \cite[Section 6]{Broue} that this
map is $\La$-linear, and so $\tr^{\La}_{\Ga}$ is a $k$-linear map
$$\Hom_{\Ga}(M,N) \xrightarrow{\tr^{\La}_{\Ga}}
\Hom_{\La}(M,N)$$ called the \emph{trace map}. Using this trace
map, the following lemma shows that the restriction map
$$\Hom_{\La}(M,N) \xrightarrow{\res^{\La}_{\Ga}} \Hom_{\Ga}(M,N)$$
is injective whenever $\mu (c_{\Ga}^{\La})$ is invertible in the
center of $\La$.

\begin{lemma}\label{restrictionhom}
Let $\La$ be a finite dimensional symmetric algebra and $\Ga$ a
parabolic subalgebra. If $\mu (c_{\Ga}^{\La})$ is invertible in
the center of $\La$, then the restriction map
$$\Hom_{\La}(M,N) \xrightarrow{\res^{\La}_{\Ga}} \Hom_{\Ga}(M,N)$$
is injective for all $\La$-modules $M$ and $N$.
\end{lemma}

\begin{proof}
Write $c_{\Ga}^{\La} = \sum_{i=1}^t x_i \otimes y_i$, and consider
the composition
$$\Hom_{\La}(M,N) \xrightarrow{\tr^{\La}_{\Ga} \circ \res^{\La}_{\Ga}} \Hom_{\La}(M,N).$$
Then
$$( \tr^{\La}_{\Ga} \circ \res^{\La}_{\Ga} )(f)(m) = \sum_{i=1}^t x_i
f(y_im) = \left ( \sum_{i=1}^t x_i y_i \right ) f(m) = \mu
(c_{\Ga}^{\La}) f(m)$$ for every $f \in \Hom_{\La}(M,N)$ and every
$m \in M$. Since $\mu (c_{\Ga}^{\La})$ is invertible, the
composition $\tr^{\La}_{\Ga} \circ \res^{\La}_{\Ga}$ is injective,
hence $\res^{\La}_{\Ga}$ is injective.
\end{proof}

Our aim in this section is to extend this result to cohomology. To do that, we need a lemma on the transitivity of trace maps. Suppose we have a parabolic chain $\La \supseteq \Gamma \supseteq \Delta$ of symmetric algebras, that is, the algebra $\Gamma$ is a parabolic subalgebra of $\La$, and $\Delta$ is a parabolic subalgebra of $\Gamma$. Since $\Gamma$ is a parabolic subalgebra of $\La$, the $\Gamma$-$\Gamma$-bimodule $\La$ is a direct sum $\La = \Gamma \oplus B_1$, where $B_1 \subseteq \Ker s$. Similarly, since $\Delta$ is a parabolic subalgebra of $\Gamma$, the $\Delta$-$\Delta$-bimodule $\Gamma$ is a direct sum $\Gamma = \Delta \oplus B_2$, where $B_2 \subseteq \Ker s$. Therefore, as a bimodule over $\Delta$, the algebra $\La$ is a direct sum $\La = \Delta \oplus B_1 \oplus B_2$, with $B_1 \oplus B_2 \subseteq \Ker s$. Consequently, the algebra $\Delta$ is a parabolic subalgebra of $\La$, and the following lemma shows the trace map is transitive in this case.

\begin{lemma}\label{transitivity}
If $\La \supseteq \Gamma \supseteq \Delta$ is a parabolic chain of symmetric algebras, then
$$\tr^{\La}_{\Ga} \circ \tr^{\Ga}_{\Delta} = \tr^{\La}_{\Delta}.$$
\end{lemma}

\begin{proof}
This follows directly from \cite[Proposition 2.11(ii)]{Linckelmann1}. With the notation used in that proposition, take $X$ and $Y$ to be the bimodules ${_{\La}\La}_{\Gamma}$ and ${_{\Ga}\Ga}_{\Delta}$, respectively.
\end{proof}

The field $k$ is obviously a parabolic subalgebra of $\La$ and
$\Ga$ when the symmetrizing form $s$ is nonzero. Moreover, the
latter is automatic if $\mu (c_{\Ga}^{\La})$ is invertible in the
center of $\La$. Namely, the relative Casimir element is by
definition the image of the unit $1_{\La}$ in $\La$ under the
composition
$$\La \xrightarrow{\phi} D( \La ) \xrightarrow{D( \mu )} D( \La
\otimes_{\Ga} \La ) \stackrel{\sim}{\longrightarrow} D( \La )
\otimes_{\Ga} D( \La ) \xrightarrow{\phi^{-1} \otimes \phi^{-1}}
\La \otimes_{\Ga} \La,$$ and the image of $1_{\La}$ under the
first map in this composition is precisely $s$. Therefore, if
$\mu (c_{\Ga}^{\La})$ is invertible, then $\La \supseteq \Ga \supseteq k$ is a parabolic chain of symmetric algebras, and the above lemma applies. Using this, we end this section with a result which extends Lemma \ref{restrictionhom} to cohomology. 

\begin{proposition}\label{restrictionext}
Let $\La$ be a finite dimensional symmetric algebra, and $\Ga$ a
parabolic subalgebra such that $\mu (c_{\Ga}^{\La})$ is invertible in $\La$.
Then for every $i$ the restriction map
$$\Ext^i_{\La}(M,N) \xrightarrow{\res^{\La}_{\Ga}} \Ext^i_{\Ga}(M,N)$$
is injective for all $\La$-modules $M$ and $N$.
\end{proposition}

\begin{proof}
Given $\La$-modules $M$ and $N$, let $f \in \Hom_{\La}(M,N)$ be a
homomorphism factoring through a projective module. Since
projective $\La$-modules are projective also as $\Ga$-modules, the
restriction $\res^{\La}_{\Ga}(f)$ factors through a projective
$\Ga$-module. Thus restriction induces a $k$-linear map
$$\uHom_{\La}(M,N) \xrightarrow{\res^{\La}_{\Ga}} \uHom_{\Ga}(M,N)$$
for stable homomorphisms. Let $g$ represent an element in
$\uHom_{\La}(M,N)$, and suppose that $\res^{\La}_{\Ga}(g)=0$ in
$\uHom_{\Ga}(M,N)$. By \cite[Lemma 3.15]{Broue}, there is a
$k$-linear map $h \in \Hom_k(M,N)$ with the property that
$\res^{\La}_{\Ga}(g) = \tr^{\Ga}_{k}(h)$. Then
$$( \tr^{\La}_{\Ga} \circ \res^{\La}_{\Ga} ) (g) = ( \tr^{\La}_{\Ga}
\circ\tr^{\Ga}_{k} ) (h) = \tr^{\La}_{k} (h)$$ by Lemma \ref{transitivity},
and so by
\cite[Lemma 3.15]{Broue} again, the map $( \tr^{\La}_{\Ga} \circ
\res^{\La}_{\Ga} ) (g)$ factors through a projective $\La$-module.
In the proof of Lemma \ref{restrictionhom} we saw that $(
\tr^{\La}_{\Ga} \circ \res^{\La}_{\Ga} ) (g) = \mu (c_{\Ga}^{\La})
g$, hence $g$ itself must factor through a projective
$\La$-module.

We have just shown that the restriction map
$$\uHom_{\La}(M,N) \xrightarrow{\res^{\La}_{\Ga}} \uHom_{\Ga}(M,N)$$
is injective for all $\La$-modules $M$ and $N$. Since the
cohomology group $\Ext^i_{\La}(M,N)$ is isomorphic to
$\uHom_{\La}( \Omega_{\La}^i(M),N)$, where $\Omega_{\La}^i(M)$
denotes the $i$th syzygy of $M$, the result follows.
\end{proof}

\section{Hecke algebras of type $A$ and symmetric groups}\label{secrepdim}

Let $\Lambda$ be a finite dimensional algebra over a field $k$,
and denote by $\mod \Lambda$ the category of finitely generated
left $\Lambda$-modules. As in the previous sections, all modules
are assumed to be finitely generated left modules. The
\emph{representation dimension} of $\Lambda$, denoted $\repdim
\Lambda$, is defined as
$$\repdim \Lambda \stackrel{\text{def}}{=} \inf \{ \gldim
\End_{\Lambda}(M) \mid M \text{ generates and cogenerates } \mod
\Lambda \},$$ where $\gldim$ denotes the global dimension of an
algebra. To say that a module generates and cogenerates $\mod
\Lambda$ means that it contains all the indecomposable projective
and injective modules as direct summands. Of course, if $\La$ is
selfinjective, these two notions coincide. It follows immediately
from the definition that a semisimple algebra is of representation
dimension zero. As mentioned, Auslander showed that the
representation dimension of a non-semisimple algebra is two if its
representation type is finite, and at least three whenever its
representation type is infinite. Thus no algebra is of
representation dimension one. Moreover, Auslander showed that the
representation dimension of a selfinjective algebra is at most its
Loewy length. Later, Iyama showed in \cite{Iyama} that the
representation dimension is finite for every finite dimensional
algebra.

The focus of this paper is the representation dimension of Hecke
algebras and group algebras of symmetric groups. Let $k$ be a
field of characteristic zero, $\ell \ge 2$ an integer and $q \in
k$ a primitive $\ell$th root of unity. Recall that the
corresponding Hecke algebra $\mathcal{H}_q(A_{n-1})$ of type $A_{n-1}$ is the $k$-algebra with generators $T_1, \dots,
T_{n-1}$ satisfying the relations
\begin{eqnarray*}
(T_i+1)(T_i-q)=0 & \text{for} & 1 \le i \le n-1 \\
T_iT_{i+1}T_i = T_{i+1}T_iT_{i+1} & \text{for} & 1 \le i \le n-2 \\
T_iT_j = T_jT_i & \text{for} & |i-j| \ge 2.
\end{eqnarray*}
If $q=1$, this is just the group algebra of the symmetric
group $S_n$, hence $\mathcal{H}_q(A_{n-1})$ is also referred to as the Hecke algebra of $S_n$.
It is well-known that a Hecke algebra is symmetric, and that its representation type depends on the number $[n/ \ell ]$ (cf.\ \cite{ErdmannNakano}). Namely, write $n= \ell m +a$, where $0 \le a \le \ell -1$ (hence $[n/ \ell ]=m$). Then $\mathcal{H}_q(A_{n-1})$ is semisimple if and only if $\ell = \infty$ or $m=0$. It is non-semisimple of finite representation type if and only if $m=1$, and of tame representation type if and only if $\ell =2$ and $n$ is either $4$ or $5$ (and then $m=2$). In all other cases, the algebra $\mathcal{H}_q(A_{n-1})$ is of wild representation type.

It was shown in \cite{BensonErdmannMikaelian} and \cite{Linckelmann} there exists a maximal $\ell$-parabolic subalgebra $\mathcal{B} \subseteq \mathcal{H}_q(A_{n-1})$ which is isomorphic to an $m$-fold tensor product
$$\mathcal{B} \simeq B_1 \otimes \cdots \otimes B_m,$$
where each $B_i$ is either semisimple or a Brauer tree algebra. Thus each $B_i$ has finite representation type. The same occurs for group algebras of certain symmetric groups. Suppose $k$ is a field of positive characteristic $p$, let $n$ be an integer with $n < p^2$, and $S_n$ the $n$th symmetric group. If we write $m = [n/p]$, then a Sylow $p$-subgroup of $S_n$ is a direct product $P = P_1 \times \cdots \times P_m$, where each $P_i$ has order $p$, and these cyclic groups have disjoint supports. The group algebra $kP$ is isomorphic to the tensor product $kP_1 \otimes \cdots \otimes kP_m$, and each $kP_i$ is isomorphic to the local algebra $k[x]/(x^p)$. In particular, each algebra $kP_i$ has finite representation type. 

Note that in both these situations, there is a canonical generator for the subalgebra. Namely, for each $1 \le i \le m$, let $M_i$ be the direct sum of a complete set of isomorphism classes of indecomposable $B_i$-modules (respectively, $kP_i$-modules). Then the $\mathcal{B}$-module (respectively, $kP$-module) $M_1 \otimes \cdots \otimes M_m$ is a generator. The following crucial lemmas show that when we induce and then restrict this module, the resulting module is contained in the additive closure of $M$. We include a proof only for the group algebra case; the proof of the Hecke algebra case is completely analogous.

\begin{lemma}\label{inducerestrictgroup}
Let $k$ be a field of positive characteristic $p$, let $n$ be an integer with $n < p^2$, and $S_n$ the $n$th symmetric group. Furthermore, let $P = P_1 \times \cdots \times P_m$ be a Sylow $p$-subgroup, where each $P_i$ has order $p$ and $m=[n/p]$. Finally,
for each $1 \le i \le m$, let $M_i$ be the direct sum of a complete set of isomorphism classes of indecomposable $kP_i$-modules, and denote the $kP$-module $M_1 \otimes \cdots \otimes M_m$ by $M$. Then the $kP$-module $kS_n \otimes_{kP} M$ belongs to $\add_{kP} M$.
\end{lemma}

\begin{proof}
By the Mackey formula, the restriction of $kS_n \otimes_{kP} M$ to $kP$ is given by
$$kS_n \otimes_{kP} M = \bigoplus_x kP \otimes_{k(P \cap xPx^{-1})} (x \otimes M),$$
where the sum is taken over a system of double coset representatives. Our module $M$ is an outer tensor product of the $kP_i$-modules $M_i$, so assume that each $P \cap xPx^{-1}$ in the formula is a direct product $P \cap xPx^{-1} = Q_1 \times \cdots \times Q_m$, with each $Q_i$ a subgroup of $P_i$. Then the $k(P \cap xPx^{-1})$-module $x \otimes M$ is again an outer tensor product of modules over the algebras $kQ_i$, hence $kP \otimes_{k(P \cap xPx^{-1})} (x \otimes M)$ is also an outer tensor product $N_1 \otimes \cdots \otimes N_m$, with each $N_i$ a module over $kP_i$. Thus, in this case, each $kP$-module $kP \otimes_{k(P \cap xPx^{-1})} (x \otimes M)$, and therefore also $kS_n \otimes_{kP} M$, belongs to $\add_{kP} M$.

We must therefore show that each $P \cap xPx^{-1}$ is a direct product $P \cap xPx^{-1} = Q_1 \times \cdots \times Q_m$, with each factor $Q_i$ a subgroup of $P_i$. Write $n=mp +a$, where $0 \le a <p$. The group $P \cap xPx^{-1}$ is a subgroup of $S_{\lambda} \cap xS_{\lambda}x^{-1}$, where $S_{\lambda}$ is a Young subgroup for the partition $\lambda = (p^m,1^a)$ of $n$, and $P \le S_{\lambda}$. The intersection $S_{\lambda} \cap xS_{\lambda}x^{-1}$ is again a Young subgroup, for a partition which refines $\lambda$. The group $P \cap xPx^{-1}$ is a $p$-subgroup for this intersection, and by Sylow's theorem it is therefore contained in a Sylow $p$-subgroup $R$ of $S_{\lambda} \cap xS_{\lambda}x^{-1}$. The group $R$ is a direct product of the Sylow subgroups of the factors, and is therefore a direct product $R= R_1 \times \cdots \times R_m$. If $R_i$ is nontrivial, then it is generated by a $p$-cycle, and distinct nontrivial $R_i,R_j$ have disjoint supports. Therefore any subgroup of $R$, in particular $P \cap xPx^{-1}$, is a direct product of groups generated by $p$-cycles, with disjoint supports.
\end{proof}

As mentioned, we do not include the proof of the Hecke algebra version of Lemma \ref{inducerestrictgroup}, since it is completely analogous to the proof just given. For more background on the relevant machinery (for instance, the Mackey formula), see \cite{DipperJames}.

\begin{lemma}\label{inducerestricthecke}
Let $\mathcal{H}_q(A_{n-1})$ be the Hecke algebra of the symmetric group
$S_n$, where $q$ is a primitive $\ell$th root of unity and the
ground field is of characteristic zero. Furthermore, let $\mathcal{B} \simeq B_1 \otimes \cdots \otimes B_m$ be a maximal $\ell$-parabolic subalgebra, where $m=[n/ \ell ]$ and each $B_i$ is either semisimple or a Brauer tree algebra. Finally, for each $i$, let $M_i$ be the direct sum of a complete set of isomorphism classes of indecomposable $B_i$-modules, and denote the $\mathcal{B}$-module $M_1 \otimes \cdots \otimes M_m$ by $M$. Then the $\mathcal{B}$-module $\mathcal{H}_q(A_{n-1}) \otimes_{\mathcal{B}} M$ belongs to $\add_{\mathcal{B}} M$.
\end{lemma}

We are now ready to prove our two main results. The first of these gives both an upper and a lower bound for the
representation dimension of a non-semisimple Hecke algebra $\mathcal{H}_q(A_{n-1})$. The main ingredient in the proof of the upper bound is Theorem \ref{gldimsubalg}.

\begin{theorem}\label{mainhecke}
Let $\mathcal{H}_q(A_{n-1})$ be the Hecke algebra of type $A_{n-1}$, where $q$ is a primitive $\ell$th root of unity and the
ground field is of characteristic zero. If $\mathcal{H}_q(A_{n-1})$ is not semisimple \emph{(}i.e.\ if $\ell$ is finite and $[n/ \ell ] \ge 1$\emph{)},
then
$$[n/ \ell ] +1 \le \repdim \mathcal{H}_q(A_{n-1}) \le 2 [n/ \ell ].$$
\end{theorem}

\begin{proof}
Let $k$ be the ground field, and write $n= \ell m +a$, where $0
\le a \le \ell -1$. We must show that
$$m+1 \le \repdim \mathcal{H}_q(A_{n-1}) \le 2m,$$
and we start with the lower bound. Throughout this proof, we
denote our Hecke algebra $\mathcal{H}_q(A_{n-1})$ by $\La$.

By \cite[Theorem 1.1]{Linckelmann}, the Hochschild cohomology ring
$$\HH^* ( \La ) = \Ext^*_{\La \otimes_k \La^{\op}}( \La, \La )$$ of $\La$ is
Noetherian, and $\Ext^*_{\La \otimes_k \La^{\op}}( \La, X )$ is a
finitely generated $\HH^* ( \La )$-module for every
$\La$-$\La$-bimodule $X$. By \cite[Proposition 2.4]{ErdmannEtAl},
the latter is equivalent to $\Ext_{\La}^*( \La / \rad \La, \La /
\rad \La )$ being a finitely generated $\HH^* ( \La )$-module.
Since the characteristic of $k$ is zero, it is a perfect field,
hence $( \La / \rad \La ) \otimes_k ( \La / \rad \La )$ is a
semisimple algebra. Therefore, by \cite[Corollary 3.6]{Bergh} (see
also \cite[Corollary 5.12]{BIKO}), the inequality
$$\dim \HH^* ( \La ) + 1 \le \repdim \La$$
holds, where $\dim \HH^* ( \La )$ is the Krull dimension of $\HH^*
( \La )$. By \cite[Theorem 1.2]{Linckelmann}, the Krull dimension
of $\HH^* ( \La )$ is $m$, hence the lower bound follows.

To prove the upper bound, let $\mathcal{B}$ be a maximal $\ell$-parabolic
subalgebra of $\La$. Then $\mathcal{B}$ is isomorphic to an
$m$-fold tensor product $\mathcal{B} \simeq B_1 \otimes \cdots
\otimes B_m,$ where each $B_i$ is either semisimple or a Brauer
tree algebra. In any case, each $B_i$ is of finite representation
type. For each $i$, let $M_i$ be the direct sum of a complete set
of isomorphism classes of indecomposable $B_i$-modules, and
consider the $\mathcal{B}$-module
$$M = M_1 \otimes \cdots \otimes M_m.$$
Then $\gldim \End_{B_i}(M_i) \le 2$, and therefore
$$\gldim \End_{\mathcal{B}} (M) = \sum_{i=1}^m \gldim
\End_{B_i}(M_i)\le 2m$$ by \cite[Corollary 3.3 and Lemma 3.4]{Xi},
since the ground field $k$ is perfect. The module $M$ is a
generator in $\mod \mathcal{B}$, and the $\mathcal{B}$-module $\La
\otimes_{\mathcal{B}} M$ belongs to $\add_{\mathcal{B}} M$ by Lemma \ref{inducerestricthecke}.
Moreover, by \cite[Theorem 2.7]{Du} the element $\mu (
c^{\La}_{\mathcal{B}} )$ is invertible in the center of $\La$, hence by Proposition \ref{restrictionext},
the restriction map
$$\Ext^i_{\La}(X,Y) \xrightarrow{\res^{\La}_{\mathcal{B}}} \Ext^i_{\mathcal{B}}(X,Y)$$
is injective for all $i$ and all $\La$-modules $X$ and $Y$.
Then by Theorem \ref{gldimsubalg} the inequality 
$$\gldim \End_{\La} ( \La \otimes_{\mathcal{B}} M ) \le \gldim \End_{\mathcal{B}} ( M
)\le 2m$$
holds. Consequently, since the $\La$-module $\La
\otimes_{\mathcal{B}} M$ is a generator by Proposition \ref{addsubalg}, 
the representation dimension of $\La$ is at most $2m$.
\end{proof}

Next, we prove the second of our main results, namely the analogue of Theorem \ref{mainhecke} for group algebras of symmetric groups. Here we use Theorem \ref{gldimseparably} when we prove the upper bound.

\begin{theorem}\label{maingroup}
Let $k$ be a perfect field of positive characteristic $p$, let $n$ be an integer with $n < p^2$, and $S_n$ the $n$th symmetric group. If $kS_n$ is not semisimple \emph{(}i.e.\ if $p \le n$\emph{)},
then
$$[n/p] +1 \le \repdim kS_n \le 2 [n/p].$$
\end{theorem}

\begin{proof}
We start with the lower bound. Since $n<p^2$, any Sylow $p$-subgroup of $S_n$ is elementary abelian of order $[n/p]$, hence this number is the $p$-rank of $S_n$. The lower bound now follows from \cite[Corollary 19]{Oppermann1}.

For the upper bound, let $P$ be a Sylow $p$-subgroup of $S_n$. Then $P = P_1 \times \cdots \times P_m$, where each $P_i$ has order $p$ and $m=[n/p]$. For each $1 \le i \le m$, let $M_i$ be the direct sum of a complete set of isomorphism classes of indecomposable $kP_i$-modules, and denote the $kP$-module $M_1 \otimes \cdots \otimes M_m$ by $M$. This module generates $\mod kP$, and as in the proof of Theorem \ref{mainhecke}, the global dimension of $\End_{kP}(M)$ is at most $2m$.

Define bimodules ${_{kS_n}X}_{kP}$ and ${_{kP}Y}_{kS_n}$ by $X= kS_n$ and $Y=kS_n$. Then the $kS_n$-$kS_n$-bimodule $kS_n$ is a direct summand of $X \otimes_{kP} Y$, so $kS_n$ separably divides $kP$. Moreover, the $kP$-module $\Hom_{kS_n}(X, X \otimes_{kP}M)$ is just $X \otimes_{kP}M$, and therefore belongs to $\add_{kP}M$ by Lemma \ref{inducerestrictgroup}. Theorem \ref{gldimseparably} now gives
$$\gldim \End_{kS_n}(X \otimes_{kP}M) \le \gldim \End_{kP}(M) \le 2m,$$
and since $X \otimes_{kP}M$ generates $\mod kS_n$ by Proposition \ref{addseparably}, the result follows.
\end{proof}

We end this section with some remarks on our main results.

\begin{remarks}
(1) Instead of Theorem \ref{gldimsubalg}, we could have used Theorem \ref{gldimseparably} to prove the upper bound in Theorem \ref{mainhecke}. Indeed, if $\mathcal{B}$ is a maximal $\ell$-parabolic subalgebra of $\mathcal{H}_q(A_{n-1})$, define bimodules ${_{\mathcal{H}_q(A_{n-1})}X}_{\mathcal{B}}$ and ${_{\mathcal{B}}Y}_{\mathcal{H}_q(A_{n-1})}$ by $X=Y= \mathcal{H}_q(A_{n-1})$. By \cite[Proposition 5.1]{Linckelmann}, the algebra $\mathcal{H}_q(A_{n-1})$ separably divides $\mathcal{B}$ through $X$ and $Y$. Thus assumption (1) in Theorem \ref{gldimseparably} holds. Assumption (2) holds by Lemma \ref{inducerestricthecke} and the arguments used at the end of the proof of Theorem \ref{maingroup}.

(2) Similarly, instead of Theorem \ref{gldimseparably}, we could have used Theorem \ref{gldimsubalg} to prove the upper bound in Theorem \ref{maingroup}. First, note that when $P$ is a Sylow $p$-subgroup of $S_n$, then $kP$ is a parabolic subalgebra of $kS_n$. When $n<p^2$, the index $|S_n|/|P|$ is not divisible by $p$, hence $\mu (c^{kS_n}_{kP})$ is invertible in $kS_n$. Therefore, by Proposition \ref{restrictionext}, assumption (2) in Theorem \ref{gldimsubalg} holds. Assumption (1) is just Lemma \ref{inducerestrictgroup}.

(3) As mentioned prior to Lemma \ref{inducerestrictgroup}, the Hecke algebra $\mathcal{H}_q(A_{n-1})$
is non-semisimple of finite representation type if and only if $[n/ \ell ]=1$, and of tame representation type if and only if $\ell =2$ and $n$ is either $4$ or $5$ (and then $[n/ \ell ]=2$). In all the other non-semisimple cases, the Hecke algebra has wild representation type. In the finite type case, Theorem \ref{mainhecke} therefore gives $2 \le \repdim \mathcal{H}_q(A_{n-1}) \le 2$, i.e.\ $\repdim \mathcal{H}_q(A_{n-1}) = 2$. This was of course to be expected: every non-semisimple finite dimensional algebra of finite representation type has representation dimension $2$. When $\mathcal{H}_q(A_{n-1})$ is tame, Theorem \ref{mainhecke} gives $3 \le \repdim \mathcal{H}_q(A_{n-1}) \le 4$, that is, the representation dimension is either $3$ or $4$. It is known that the representation dimension is $3$ in this case. Namely, by \cite{SchrollSnashall}, there are two Morita equivalence classes of blocks, represented by $\mathcal{H}_q(A_3)$ and $\mathcal{H}_q(A_4)$. It is shown in \cite{ErdmannNakano} that these algebras are special biserial, and so by \cite[Corollary 1.3]{ErdmannHolmIyamaSchroer} they are both of representation dimension $3$.

(4) The group algebra $kS_n$ is non-semisimple when $p \le n$, and when $n<p^2$ it has finite representation type precisely when $[n/p]=1$. For in this case, the Sylow $p$-subgroups of $S_n$ are cyclic, and by a classical result of Higman this is equivalent to $kS_n$ having finite type. As was the case for Hecke algebras, Theorem \ref{maingroup} gives $\repdim kS_n =2$ in this case.
\end{remarks}

\section{Hecke algebras of types $B$ and $D$}

In this final section, we consider Hecke algebras 
of types $B_n$ ($n\ge 2$) and $D_n$ ($n\ge 4$).
These are associated to Coxeter groups 
$W(B_n)=\cS_2\wr \cS_n$ and $W(D_n)=W(B_n)\cap \Alt_{2n}$, where $\Alt_{2n}$
is the alternating group of degree $2n$.
As before, the ground field is assumed to be of characteristic zero.

For type $B_n$, the Hecke algebra involves two parameters $q$ and $Q$,
and we write $\cH_{Q,q}(B_n)$. This is the $k$-algebra with generators $T_0, T_1, \dots,
T_{n-1}$ satisfying the relations
\begin{eqnarray*}
(T_0+1)(T_0-Q)=0 \\
T_0T_1T_0T_1 = T_1T_0T_1T_0 \\
(T_i+1)(T_i-q)=0 & \text{for} & 1 \le i \le n-1 \\
T_iT_{i+1}T_i = T_{i+1}T_iT_{i+1} & \text{for} & 1 \le i \le n-2 \\
T_iT_j = T_jT_i & \text{for} & |i-j| \ge 2.
\end{eqnarray*}
Note that there is a natural inclusion of the Hecke algebra $\cH_q(A_{n-1})$ of type $A_{n-1}$ into $\cH_{Q,q}(B_n)$, taking the generator $T_i$ of $\cH_q(n)$ (for $1 \le i \le n-1$) to the generator $T_i$ of $\cH_{Q,q}(B_n)$,
ignoring the element $T_0$. The algebra $\cH_{Q,q}(B_n)$ is free as a left/right module over $\cH_q(A_{n-1})$. 

For Hecke algebras of type $D_n$, there is just one parameter $q$, 
and we write $\cH_q(D_n)$. This is the $k$-algebra with generators $T_0, T_1, \dots,
T_{n-1}$ satisfying the relations
\begin{eqnarray*}
T_0T_2T_0 = T_2T_0T_2 \\
T_0T_i = T_iT_0 & \text{for} & 1 \le i \le n-1, i \neq 2 \\
(T_i+1)(T_i-q)=0 & \text{for} & 0 \le i \le n-1 \\
T_iT_{i+1}T_i = T_{i+1}T_iT_{i+1} & \text{for} & 1 \le i \le n-2 \\
T_iT_j = T_jT_i & \text{for} & 1 \le i \le j-2 \le n-3.
\end{eqnarray*}
As with Hecke algebras of type $B$, there is a natural inclusion of $\cH_q(A_{n-1})$ into $\cH_q(D_n)$: 
the generator $T_i$ of $\cH_q(n)$ (for $1 \le i \le n-1$) maps to the generator $T_i$ of $\cH_q(D_n)$,
ignoring the element $T_0$. Moreover, the algebra $\cH_q(D_n)$ is free as a left/right module over $\cH_q(A_{n-1})$. 

We shall first establish lower bounds for the representation dimensions of Hecke algebras of types $B$ and $D$. To do this, we compute lower bounds for the dimensions of the stable module categories involved, viewed as triangulated categories. Namely, if $\cH$ is any Hecke algebra, then $\cH$ is symmetric, in particular selfinjective. The stable module category $\stmod \cH$ is then triangulated, with suspension functor $\Omega_{\cH}^{-1} \colon \stmod \cH \to \stmod \cH$. Now let $\mathcal{T}$ be an arbitrary triangulated category with suspension functor $\Sigma \colon \mathcal{T} \to \mathcal{T}$, and $\mathcal{C}$ and $\mathcal{D}$ subcategories of $\mathcal{T}$. We denote by $\thick^1_{\mathcal{T}} ( \mathcal{C} )$ the full subcategory of $\mathcal{T}$ consisting of all the direct summands of finite direct sums of suspensions of objects in $\mathcal{C}$. Furthermore, we denote by $\mathcal{C} \ast \mathcal{D}$ the full subcategory of $\mathcal{T}$
consisting of objects $M$ such that there exists a distinguished triangle
$$C \to M \to D \to \Sigma C$$
in $\mathcal{T}$, with $C \in \mathcal{C}$ and $D \in \mathcal{D}$. Now for each $n \ge 2$,
define inductively $\thick^n_{\mathcal{T}} ( \mathcal{C} )$ by
$$\thick^n_{\mathcal{T}} ( \mathcal{C} ) \stackrel{\text{def}}{=} \thick_{\mathcal{T}}^1 \left
( \thick^{n-1}_{\mathcal{T}} ( \mathcal{C} ) \ast \thick^1_{\mathcal{T}} ( \mathcal{C} ) \right ).$$
Then the \emph{dimension} of $\mathcal{T}$, denoted $\dim \mathcal{T}$, is defined as
$$\dim \mathcal{T} \stackrel{\text{def}}{=} \inf \{ n \ge 0 \mid \exists \text{ an object } G \in \mathcal{T} \text{ such that } \mathcal{T} = \thick^{n+1}_{\mathcal{T}} (G) \}.$$
This notion was introduced by Rouquier in \cite{Rouquier1}, precisely in order to establish lower bounds for the representation dimension of certain algebras, namely exterior algebras. 

\begin{lemma}\label{rouquier}\cite[Proposition 3.9]{Rouquier1}
If $\La$ is a finite dimensional non-semisimple selfinjective algebra, then $\dim \left ( \stmod \La \right ) +2 \le \repdim \La$.
\end{lemma}

Thus, in order to compute lower bounds for the representation dimensions Hecke algebras of types $B$ and $D$, we compute lower bounds for the dimensions of their stable module categories. We include also Hecke algebras of type $A$.

\begin{lemma}\label{dimension}
Let $k$ be a field of characteristic zero, and $q \in k$ a primitive $\ell$th root of unity, where $\ell$ is finite. Furthermore, let $n$ be an integer, and $\cH$ either a Hecke algebra $\cH_q(A_{n-1})$ of type $A_{n-1}$, a Hecke algebra $\cH_{Q,q}(B_n)$ of type $B_n$, or a Hecke algebra $\cH_q(D_n)$ of type $D_n$.
Then $\dim \left ( \stmod \cH \right ) \ge [n/ \ell ]-1$.
\end{lemma}

\begin{proof}
The Hecke algebra $\cH_{q}(A_{n-1})$ of type $A_{n-1}$ is a subalgebra of $\cH$, but also a factor algebra in a natural way (factor out the generator $T_0$). Therefore, there is a diagram 
$$\cH_{q}(A_{n-1}) \xrightarrow{i} \cH \xrightarrow{\pi} \cH_{q}(A_{n-1})$$
of $k$-algebra homomorphisms, in which the composition $\pi \circ i$ is the identity on $\cH_{q}(A_{n-1})$. Since $\cH$ is projective as a left $\cH_{q}(A_{n-1})$-module, the inequality
$$\dim \left ( \stmod \cH_{q}(A_{n-1}) \right ) \le \dim \left ( \stmod \cH \right )$$
holds by \cite[Lemma 2.3]{BerghOppermann1}. It is therefore enough to prove the result for $\cH_{q}(A_{n-1})$.

Denote our algebra $\cH_{q}(A_{n-1})$ by $\La$. Since the ground field $k$ is of characteristic zero, it is a perfect field, hence the algebra $( \La / \rad \La ) \otimes_k ( \La / \rad \La )$ is semisimple. As we saw in the proof of Theorem \ref{mainhecke}, the Hochschild cohomology ring $\HH^*( \La )$ is Noetherian of Krull dimension $[n/ \ell ]$, and $\Ext_{\La}^*( \La / \rad \La, \La / \rad \La )$ is a finitely generated $\HH^*( \La )$-module. Moreover, in the proof of \cite[Corollary 3.6]{Bergh}, it was shown that the Krull dimension of the Hochschild cohomology ring equals the complexity of the $\La$-module $\La / \rad \La$. It then follows from \cite[Theorem 3.2]{Bergh} that $\dim \left ( \stmod \La \right ) \ge [n/ \ell ]-1$.
\end{proof}

Combining this lemma with Lemma \ref{rouquier}, we obtain the lower bounds for the representation dimensions of Hecke algebras of types $B$ and $D$. The same bound holds for Hecke algebras of types $A_{n-1}$, $B_n$ and $D_n$.

\begin{theorem}\label{lowerBD}
Let $k$ be a field of characteristic zero, and $q \in k$ a primitive $\ell$th root of unity, where $\ell$ is finite. Furthermore, let $n$ be an integer such that $[n/ \ell ] \ge 1$, and $\cH$ either a Hecke algebra $\cH_{Q,q}(B_n)$ of type $B_n$, or a Hecke algebra $\cH_q(D_n)$ of type $D_n$. Then $\repdim \cH \ge [n/ \ell ]+1$.
\end{theorem}

\begin{remark}
There is a common generalization of Hecke algebras of types $A$ and $B$, namely the \emph{Ariki-Koike} algebras. Let $q, Q_1, \dots, Q_t$ be elements of $k$, and denote the sequence $( Q_1, \dots, Q_t )$ by ${\bf{Q}}$. The corresponding Ariki-Koike algebra $\cH_{{\bf{Q}},q}(n)$ is the $k$-algebra with generators $T_0, T_1, \dots,
T_{n-1}$ satisfying the relations
\begin{eqnarray*}
\prod_{i=1}^t(T_0-Q_i)=0 \\
T_0T_1T_0T_1 = T_1T_0T_1T_0 \\
(T_i+1)(T_i-q)=0 & \text{for} & 1 \le i \le n-1 \\
T_iT_{i+1}T_i = T_{i+1}T_iT_{i+1} & \text{for} & 1 \le i \le n-2 \\
T_iT_j = T_jT_i & \text{for} & |i-j| \ge 2.
\end{eqnarray*}
When ${\bf{Q}} = (1)$, this is just the Hecke algebra $\cH_{q}(A_{n-1})$ of type $A_{n-1}$, whereas when ${\bf{Q}} = (-1,Q)$ we retrieve the Hecke algebra $\cH_{Q,q}(B_n)$ of type $B_n$. As with Hecke algebras of types $B_n$ and $D_n$, the algebra $\cH_{q}(A_{n-1})$ is a subalgebra of $\cH_{{\bf{Q}},q}(n)$, and the latter is free over $\cH_{q}(A_{n-1})$. Moreover, there is a diagram 
$$\cH_{q}(A_{n-1}) \xrightarrow{i} \cH_{{\bf{Q}},q}(n) \xrightarrow{\pi} \cH_{q}(A_{n-1})$$
of $k$-algebra homomorphisms, in which the composition $\pi \circ i$ is the identity on $\cH_{q}(A_{n-1})$. As in the proof of Lemma \ref{dimension}, we obtain the inequalities 
$$\dim \left ( \stmod \cH_{{\bf{Q}},q}(n) \right ) \ge \dim \left ( \stmod \cH_{q}(A_{n-1}) \right ) \ge [n/ \ell ]-1,$$
and so the representation dimension of $\cH_{{\bf{Q}},q}(n)$ is also bounded below by $[n/ \ell ]+1$.
\end{remark}

Finally, we turn to upper bounds for the representation dimensions of Hecke algebras of types $B$ and $D$. The situation is more complicated than for type $A$, and our method depends on whether certain polynomial expressions in the parameters are nonzero. First we treat type $B_n$ for $n\ge 2$: we set  
$$f_n(Q,q)= \prod_{i=1-n}^{n-1}(Q+q^i).$$ 
By a result of Dipper and James, when $f_n(Q,q)$ is nonzero, then the algebra $\cH_{Q,q}(B_n)$ is Morita equivalent to a product of tensor products of Hecke algebras of type $A$. Of course, the condition that $f_n(Q,q)$ be nonzero is equivalent to the condition $Q+q^i \neq 0$ for $1-n \le i \le n-1$.

\begin{lemma}\cite[Theorem 4.17]{DipperJames2}\label{moritaB}
If $f_n(Q,q)$ is nonzero, then the Hecke algebra
$\cH_{Q,q}(B_n)$ is Morita equivalent to the algebra
\[ \prod _{j=0}^n \cH_q(A_{j-1})\otimes_k\cH_q(A_{n-j-1}). \]
\end{lemma}

Using this result, we obtain the upper bound for the representation dimension in type $B$, provided the relevant polynomial is invertible.

\begin{theorem}\label{upperB}
Let $k$ be a field of characteristic zero, and $q \in k$ a primitive $\ell$th root of unity, where $\ell$ is finite. Furthermore, let $n$ be an integer, and $\cH_{Q,q}(B_n)$ a Hecke algebra of type $B_n$. If $f_n(Q,q)$ is nonzero, then $\repdim \cH_{Q,q}(B_n) \le 2[n/ \ell ]$.
\end{theorem}

\begin{proof}
In general, the representation dimension of a direct product of algebras equals the maximum of the representation dimensions of the factors. Therefore, by Lemma \ref{moritaB}, the representation dimension of $\cH_{Q,q}(B_n)$ equals 
$$\max \{ \repdim \left ( \cH_q(A_{j-1})\otimes_k\cH_q(A_{n-j-1}) \right ) \mid 0 \le j \le n \}.$$
Now by \cite[Theorem 3.5]{Xi}, the representation dimension of $\cH_q(A_{j-1})\otimes_k\cH_q(A_{n-j-1})$ is at most 
$$\repdim \cH_q(A_{j-1}) + \repdim \cH_q(A_{n-j-1}),$$
and by Theorem \ref{mainhecke} this sum is at most $2( [j/ \ell ] + [(n-j)/ \ell ])$ (of course, the upper bound in Theorem \ref{mainhecke} holds without the assumption that $[n/ \ell ] \ge 1$: if $[n/ \ell ] =0$ then the Hecke algebra is semisimple). Since $[j/ \ell ] + [(n-j)/ \ell ] \le [n/ \ell ]$, we are done.
\end{proof}

Next, we consider type $D_n$ for $n\ge 4$, and here we set 
$$g_n(q)= 2\prod_{i=1} ^{n-1}(1+q^i).$$
The following theorem, which is analogous to the one for type $B$, is implicit in \cite[Theorems 3.6 and 3.7]{Pallikaros}, and is made explicit in \cite{Hu}.

\begin{lemma}\label{moritaD}
If $g_n(q)$ is nonzero and $n$ is odd,
then $\cH_q(D_n)$ is Morita equivalent to the algebra
\[ \prod_{j=(n+1)/2}^n \cH_q(A_{j-1})\otimes_k \cH_q(A_{n-j-1}). \]
\end{lemma}

There is a  corresponding theorem when $n$ is even, a result which was proved in \cite{Hu}.
However, the Morita description in this case has a direct factor which is not a tensor product of type A Hecke algebras, and therefore we cannot deal with the case when $n$ is even. But when $n$ is odd we obtain an upper bound. We skip the proof since it is exactly the same as that of Theorem \ref{upperB}.

\begin{theorem}\label{upperD}
Let $k$ be a field of characteristic zero, and $q \in k$ a primitive $\ell$th root of unity, where $\ell$ is finite. Furthermore, let $n$ be an odd integer, and $\cH_q(D_n)$ a Hecke algebra of type $D_n$. If $g_n(q)$ is nonzero, then $\repdim \cH_q(D_n) \le 2[n/ \ell ]$.
\end{theorem}

\end{document}